\documentclass[12pt]{article}
\usepackage{amsmath,amsfonts,latexsym,amsthm,amssymb}
\numberwithin{equation}{section}

\DeclareMathOperator{\const}{const}
\DeclareMathOperator{\R}{Re}
\begin{document}

\newtheorem{lem}{Lemma}
\newtheorem{teo}{Theorem}
\newtheorem*{defin}{Definition}
\newtheorem*{exam}{Example}
\pagestyle{plain}
\title{Radial Solutions of Non-Archimedean Pseudo-Differential Equations}
\author{Anatoly N. Kochubei\footnote{This work was supported in part by Grant No. 01-01-12 of the National Academy of Sciences of Ukraine under the program of joint Ukrainian-Russian projects.}
\\ \footnotesize Institute of Mathematics,\\
\footnotesize National Academy of Sciences of Ukraine,\\
\footnotesize Tereshchenkivska 3, Kiev, 01601 Ukraine\\
\footnotesize E-mail: \ kochubei@i.com.ua}
\date{}
\maketitle

\vspace*{3cm}
\begin{abstract}
We consider a class of equations with the fractional differentiation operator $D^\alpha$, $\alpha >0$, for complex-valued functions $x\mapsto f(|x|_K)$ on a non-Archimedean local field $K$ depending only on the absolute value $|\cdot |_K$. We introduce a right inverse $I^\alpha$ to $D^\alpha$, such that the change of an unknown function $u=I^\alpha v$ reduces the Cauchy problem for an equation with $D^\alpha$ (for radial functions) to an integral equation whose properties resemble those of classical Volterra equations. This contrasts much more complicated behavior of $D^\alpha$ on other classes of functions.
\end{abstract}
\vspace{2cm}


\newpage
\section{Introduction}

Pseudo-differential equations for complex-valued functions defined on a non-Archimedean local field are among the central objects of contemporary harmonic analysis and mathematical physics; see the monographs \cite{VVZ}, \cite{K}, \cite{AKS}, and the survey \cite{Z}.

The simplest example is the fractional differentiation operator $D^\alpha$, $\alpha >0$, on the field $\mathbb Q_p$ of $p$-adic numbers (here $p$ is a prime number). It can be defined as a pseudo-differential operator with the symbol $|\xi |_p^\alpha$ where $|\cdot |_p$ is the $p$-adic absolute value or, equivalently, as an appropriate convolution operator.

Already in this case, as it was first shown by Vladimirov (see \cite{VVZ}), properties of the $p$-adic pseudo-differential operator are much more complicated than those of its classical counterpart. It suffices to say that, as an operator on $L_2(\mathbb Q_p)$, it has a point spectrum of infinite multiplicity. Considering a simple ``formal'' evolution equation with the operator $D^\alpha$ in the $p$-adic time variable $t$, Vladimirov (see \cite{V03}) noticed that such an equation does not possess a fundamental solution.

At the same time, it was found in \cite{K08} that some of the evolution equations of the above kind behave reasonably, if one considers only solutions depending on $|t|_p$. This observation has led to the concept of a non-Archimedean wave equation possessing various properties resembling those of classical hyperbolic equations, up to the Huygens principle.

In this paper we consider the Cauchy problem for a class of equations like
\begin{equation}
D^\alpha u+a(|x|_p)u=f(|x|_p),\quad x \in \mathbb Q_p,
\end{equation}
assuming that a solution is looked for in the class of radial functions, $u=u(|x|_p)$; the accurate definition of $D^\alpha$ and assumptions on $a,f$ are given below. This Cauchy problem is reduced to an integral equation resembling classical Volterra equations. It appears that the equation (1.1) and its generalizations considered on radial functions constitute $p$-adic counterparts of ordinary differential equations.

The author is grateful to M. M. Malamud for a helpful discussion.

\section{Preliminaries}

{\bf 2.1. Local fields.} Let $K$ be a non-Archimedean local field,
that is a non-discrete totally disconnected locally compact
topological field. It is well known that $K$ is isomorphic either
to a finite extension of the field $\mathbb Q_p$ of $p$-adic
numbers (if $K$ has characteristic 0), or to the field of formal
Laurent series with coefficients from a finite field, if $K$ has
a positive characteristic. For a summary of main notions and results
regarding local fields see, for example, \cite{K}.

Any local field $K$ is endowed with an absolute value $|\cdot |_K$,
such that $|x|_K=0$ if and only if $x=0$, $|xy|_K=|x|_K\cdot |y|_K$,
$|x+y|_K\le \max (|x|_K,|y|_K)$. Denote $O=\{ x\in K:\ |x|_K\le 1\}$,
$P=\{ x\in K:\ |x|_K<1\}$. $O$ is a subring of $K$, and $P$ is an ideal
in $O$ containing such an element $\beta$ that
$P=\beta O$. The quotient ring $O/P$ is
actually a finite field; denote by $q$ its cardinality. We will
always assume that the absolute value is
normalized, that is $|\beta |_K=q^{-1}$. The normalized absolute
value takes the values $q^N$, $N\in \mathbb Z$. Note that for $K=\mathbb Q_p$
we have $\beta =p$ and $q=p$; the $p$-adic absolute value is normalized.

Denote by $S\subset O$ a
complete system of representatives of the residue classes from
$O/P$. Any nonzero element $x\in K$ admits the canonical
representation in the form of the convergent series
\begin{equation}
x=\beta^{-n}\left( x_0+x_1\beta +x_2\beta^2+\cdots \right)
\end{equation}
where $n\in \mathbb Z$, $|x|_K=q^n$, $x_j\in S$, $x_0\notin P$. For
$K=\mathbb Q_p$, one may take $S=\{ 0,1,\ldots ,p-1\}$.

The additive group of any local field is self-dual, that is if
$\chi$ is a fixed non-constant complex-valued additive character of
$K$, then any other additive character can be written as
$\chi_a(x)=\chi (ax)$, $x\in K$, for some $a\in K$. Below we assume that $\chi$ is a rank
zero character, that is $\chi (x)\equiv 1$ for $x\in O$, while
there exists such an element $x_0\in K$ that $|x_0|_K=q$ and $\chi
(x_0)\ne 1$.

The above duality is used in the definition of the Fourier
transform over $K$. Denoting by $dx$ the Haar measure on the
additive group of $K$ (normalized in such a way that the measure
of $O$ equals 1) we write
$$
\widetilde{f}(\xi )=\int\limits_K\chi (x\xi )f(x)\,dx,\quad \xi
\in K,
$$
where $f$ is a complex-valued function from $L_1(K)$. As usual, the Fourier
transform $\mathcal F$ can be extended from $L_1(K)\cap L_2(K)$ to a
unitary operator on $L_2(K)$. If $\mathcal F f=\widetilde{f}\in L_1(K)$, we
have the inversion formula
$$
f(x)=\int\limits_K\chi (-x\xi )\widetilde{f}(\xi )\,d\xi .
$$

\bigskip
{\bf 2.2. Integration formulas.} As in the real analysis, there are many well-known formulas for integrals of complex-valued functions defined on subsets of a local field. There exist even tables of such integrals \cite{V03}. Note that formulas for integrals on $\mathbb Q_p$ and its subsets, as a rule, carry over to the general case, if one substitutes the normalized absolute value for $|\cdot |_p$ and $q$ for $p$.

Here we collect some formulas used in this work.

\begin{equation}
\int\limits_{|x|_K\le q^n}|x|_K^{\alpha -1}\,dx=\frac{1-q^{-1}}{1-q^{-\alpha }}q^{\alpha n};\quad \text{here and below $n\in \mathbb Z,\alpha >0$}.
\end{equation}
\begin{equation}
\int\limits_{|x|_K=q^n}|x-a|_K^{\alpha -1}\,dx=\frac{q-2+q^{-\alpha}}{q(1-q^{-\alpha})}|a|_K^\alpha ,\quad |a|_K=q^n.
\end{equation}
\begin{equation}
\int\limits_{|x|_K\le q^n}\log |x|_K\,dx=\left( n-\frac1{q-1}\right) q^n\log q.
\end{equation}
\begin{equation}
\int\limits_{|x|_K=q^n}\log |x-a|_K\,dx=\left[ \left( 1-\frac1q \right)\log |a|_K-\frac{\log q}{q-1}\right] |a|_K,\quad |a|_K=q^n.
\end{equation}
\begin{equation}
\int\limits_{|x|_K\le q^n}dx=q^n;\quad \int\limits_{|x|_K=q^n}dx=\left( 1-\frac1q \right)q^n.
\end{equation}
\begin{equation}
\int\limits_{|x|_K=q^n,x_0=k_0}dx=q^{n-1},\quad 0\ne k_0\in S
\end{equation}
(the restriction $x_0=k_0$ is in the sense of the canonical representation (2.1)).
\begin{equation}
\int\limits_{|x|_K=q^n,x_0\ne k_0}dx=\left( 1-\frac2q \right)q^n.
\end{equation}

\bigskip
{\bf 2.3. Test functions and distributions.} A function $f:\ K\to \mathbb C$ is said to be locally constant, if there exists such an integer $l$ that for any $x\in K$
$$
f(x+x')=f(x), \quad \text{whenever $|x'|\le q^{-l}$}.
$$
The smallest number $l$ with this property is called the exponent
of local constancy of the function $f$.

Let $\mathcal D(K)$ be the set of all locally constant
functions with compact supports; it is a vector space
over $\mathbb C$ with the topology of double inductive limit
$$
\mathcal D(K)=\varinjlim_{N\to \infty }\varinjlim_{l\to \infty }\mathcal D_N^l
$$
where $\mathcal D_N^l$ is the finite-dimensional space of functions supported in the ball
$B_N=\big\{ x\in K:$ $|x|\le q^N\big\}$ and having the exponents of local constancy $\le l$.
The strong conjugate space $\mathcal D'(K)$ is called the space of
Bruhat-Schwartz distributions.

The Fourier transform preserves the space $\mathcal D(K)$. Therefore the Fourier transform of a distribution defined
by duality acts continuously on $\mathcal D'(K)$. As in the case of $\mathbb R^n$, there exists a
well-developed theory of distributions over local fields; it includes
such topics as convolution, direct product, homogeneous distributions etc (see \cite{VVZ}, \cite{K}, \cite{AKS}). In connection with homogeneous distributions, it is useful to introduce the subspaces of $\mathcal D(K)$:
$$
\Psi (K)=\left\{ \psi \in \mathcal D(K):\ \psi (0)=0\right\} ,
$$
$$
\Phi (K)=\left\{ \varphi \in \mathcal D(K):\ \int\limits_K
\varphi (x)\,dx=0\right\} .
$$

The Fourier transform $\mathcal F$ is a linear isomorphism from $\Psi (K)$ onto $\Phi (K)$,
thus also from $\Phi'(K)$ onto $\Psi'(K)$. The spaces $\Phi (K)$ and $\Phi'(K)$ are called the Lizorkin spaces (of the second kind) of test functions and distributions respectively; see \cite{AKS}. Note that two distributions differing by a constant summand coincide as elements of $\Phi'(K)$.

\section{Fractional Differentiation and Integration Operators}

{\bf 3.1. Riesz kernels and operators generated by them.} On a test function $\varphi \in \mathcal D(K)$, the fractional differentiation operator $D^\alpha$, $\alpha >0$, is defined as
\begin{equation}
\left( D^\alpha \varphi \right) (x)=\mathcal F^{-1}\left[ |\xi |^\alpha
(\mathcal F (\varphi ))(\xi )\right] (x).
\end{equation}
However $D^\alpha$ does not act on the space $\mathcal D(K)$,
since the function $\xi \mapsto |\xi |^\alpha$ is not locally
constant. On the other hand, $D^\alpha :\ \Phi (K)\to \Phi (K)$
and  $D^\alpha :\ \Phi'(K)\to \Phi'(K)$; see \cite{AKS}, and that was a motivation to introduce these spaces.

The operator $D^\alpha$ can also be represented as a hypersingular integral operator:
\begin{equation}
\left( D^\alpha \varphi \right) (x)=\frac{1-q^\alpha }{1-q^{-\alpha
-1}}\int\limits_K |y|^{-\alpha -1}[\varphi (x-y)-\varphi (x)]\,dy
\end{equation}
\cite{VVZ,K}.
In contrast to (3.1), the expression in the right of (3.2) makes sense for wider classes of functions. Below we study this in detail for the case of radial functions.

The expression in (3.2) is in fact a convolution $f_{-\alpha}*\varphi$ where $f_s$, $s\in \mathbb C$, $s\not \equiv 1\pmod{\dfrac{2\pi i}{\log q}\mathbb Z}$, is the Riesz kernel defined first for $\R s>0$, $s\ne 1$, as
$$
f_s(x)=\frac{|x|_K^{s-1}}{\Gamma_K(s)},\quad \Gamma_K(s)=\frac{1-q^{s-1}}{1-q^{-s}},
$$
and then extended meromorphically, as a distribution from $\mathcal D'(K)$ given by
$$
\langle f_s,\varphi \rangle =\frac{1-q^{-1}}{1-q^{s-1}}\varphi (0)+\frac{1-q^{-s}}{1-q^{s-1}}\left[ \int\limits_{|x|_K>1}\varphi (x)\frac{dx}{|x|_K^{1-s}}+\int\limits_{|x|_K\le 1}(\varphi (x)-\varphi (0))\frac{dx}{|x|_K^{1-s}}\right],
$$
$s\ne 0,s\not \equiv 1\pmod{\dfrac{2\pi i}{\log q}\mathbb Z}$.
For $s=0$, $f_0(x)=\delta (x)$. For $s\equiv 1\pmod{\dfrac{2\pi i}{\log q}\mathbb Z}$, we define
$$
f_s(x)=\frac{1-q}{\log q}\log |x|_K.
$$

It is well known that $f_s*f_t=f_{s+t}$ in the sense of distributions from $\mathcal D'(K)$, if $s,t,s+t\not \equiv 1\pmod{\dfrac{2\pi i}{\log q}\mathbb Z}$. If these kernels are considered as distributions from $\Phi'(K)$, then $f_s*f_t=f_{s+t}$ for all $s,t\in \mathbb C$ \cite{AKS}. In view of this identity, it is natural to define the operator $D^{-\alpha}$, $\alpha >0$, setting
\begin{equation}
\left( D^{-\alpha} \varphi \right) (x)=(f_\alpha *\varphi )(x)=\frac{1-q^{-\alpha} }{1-q^{\alpha
-1}}\int\limits_K |x-y|_K^{\alpha -1}\varphi (y)\,dy,\quad \varphi \in \mathcal D(K),\ \alpha \ne 1,
\end{equation}
and
\begin{equation}
\left( D^{-1}\varphi \right) (x)=\frac{1-q}{q\log q}\int\limits_K \log |x-y|_K\varphi (y)\,dy.
\end{equation}
Then $D^\alpha D^{-\alpha}=I$ on $\mathcal D(K)$, if $\alpha \ne 1$. This property remains valid on $\Phi (K)$ also for $\alpha =1$.

The above notions and results are well known; see \cite{VVZ}, \linebreak \cite{AKS}. We come to new phenomena considering the case of radial functions.

\bigskip
{\bf 3.2. Operators on radial functions.} Let $u$ be a radial function, that is $u=u(|x|_K)$, $x\in K$. Let us find an explicit expression of $D^\alpha u$, $\alpha >0$. Below we write $d_\alpha =\dfrac{1-q^\alpha}{1-q^{-\alpha -1}}$. For $x\in K$, we denote by $x_0$ the element from $S\subset O$ appearing in the representation (2.1).

\medskip
\begin{lem}
If a function $u=u(|x|_K)$ is such that
\begin{equation}
\sum\limits_{k=-\infty}^m q^k\left| u(q^k)\right| <\infty ,\quad \sum\limits_{l=m}^\infty q^{-\alpha l}\left| u(q^l)\right| <\infty,
\end{equation}
for some $m\in \mathbb Z$, then for each $n\in \mathbb Z$ the expression in the right-hand side of (3.2) with $\varphi (x)=u(|x|_K)$ exists for $|x|_K=q^n$, depends only on $|x|_K$, and
\begin{multline}
(D^\alpha u)(q^n)=d_\alpha \left(1-\frac1q \right)q^{-(\alpha +1)n}\sum\limits_{k=-\infty}^{n-1} q^ku(q^k) +q^{-\alpha n-1}\frac{q^\alpha +q-2}{1-q^{-\alpha -1}}u(q^n)\\
+d_\alpha \left(1-\frac1q \right)\sum\limits_{l=n+1}^\infty q^{-\alpha l}u(q^l).
\end{multline}
\end{lem}

\medskip
{\it Proof}. We find, using the ultrametric properties of the absolute value, that
$$
(D^\alpha u)(x)=d_\alpha \int\limits_{|y|_K\ge |x|_K}|y|_K^{-\alpha -1}\left[ u(|x-y|_K)-u(|x|_K)\right]\,dy.
$$
If $|y|_K=|x|_K$ and $y_0\ne x_0$, then the integrand vanishes. Therefore by (2.6),
\begin{multline*}
(D^\alpha u)(x)=d_\alpha \sum\limits_{k=-\infty}^{n-1}\int\limits_{|y-x|_K=q^k}|x|_K^{-\alpha -1}\left[ u(q^k)-u(q^n)\right]\,dy\\
+d_\alpha \sum\limits_{l=n+1}^\infty \int\limits_{|y|_K=q^l}q^{-l(\alpha +1)}\left[ u(q^k)-u(q^n)\right]\,dy\\
=d_\alpha \left(1-\frac1q \right)q^{-(\alpha +1)n}\sum\limits_{k=-\infty}^{n-1} q^k\left[ u(q^k)-u(q^n)\right]
+d_\alpha \left(1-\frac1q \right)\sum\limits_{l=n+1}^\infty q^{-\alpha l}\left[ u(q^l)-u(q^n)\right].
\end{multline*}

It is clear from this expression that $(D^\alpha u)(x)$, $|x|_K=q^n$, depends only on $|x|_K$. After elementary transformations we get (3.6). $\qquad \blacksquare$

\medskip
\begin{defin}
We say that the action $D^\alpha u$, $\alpha >0$, on a radial function $u$ is defined in the strong sense, if the function $u$ satisfies (3.5), so that $D^\alpha u(|x|_K)$, $|x|_K\ne 0$, is given by (3.6), and there exists the limit
$$
D^\alpha u(0)\stackrel{\text{def}}{=}\lim\limits_{x\to 0}D^\alpha u(|x|_K).
$$
\end{defin}

\medskip
It is evident from (3.2) that $D^\alpha$ annihilate constant functions (recall that in $\Phi'(K)$ they are equivalent to zero). Therefore $D^{-\alpha }$ is not the only possible choice of the right inverse to $D^\alpha$. In particular, we will use
\begin{equation}
(I^\alpha \varphi )(x)=(D^{-\alpha}\varphi )(x)-(D^{-\alpha}\varphi )(0).
\end{equation}
This is defined initially for $\varphi \in \mathcal D(K)$. It is seen from (3.3), (3.4), and the ultrametric property of the absolute value that
\begin{equation}
(I^\alpha \varphi )(x)=\frac{1-q^{-\alpha}}{1-q^{\alpha -1}}\int\limits_{|y|_K\le |x|_K}\left( |x-y|_K^{\alpha -1}-|y|_K^{\alpha -1}\right) \varphi (y)\,dy,\quad \alpha \ne 1,
\end{equation}
and
\begin{equation}
(I^1\varphi )(x)=\frac{1-q}{q\log q}\int\limits_{|y|_K\le |x|_K}\left( \log |x-y|_K-\log |y|_K\right) \varphi (y)\,dy.
\end{equation}
In contrast to (3.3) and (3.4), in (3.8) and (3.9) the integrals are taken, for each fixed $x\in K$, over bounded sets.

Let us calculate $I^\alpha u$ for a radial function $u=u(|x|_K)$. Obviously, $(I^\alpha u)(0)=0$ whenever $I^\alpha$ is defined.

\medskip
\begin{lem}
Suppose that
$$
\sum\limits_{k=-\infty}^m \max \left( q^k,q^{\alpha k}\right) \left| u(q^k)\right| <\infty ,\quad \text{if $\alpha \ne 1$},
$$
and
$$
\sum\limits_{k=-\infty}^m |k|q^k \left| u(q^k)\right| <\infty ,\quad \text{if $\alpha =1$},
$$
for some $m\in \mathbb Z$. Then $I^\alpha u$ exists, it is a radial function, and for any $x\ne 0$,
\begin{equation}
(I^\alpha u)(|x|_K)=q^{-\alpha}|x|_K^\alpha u(|x|_K)+\frac{1-q^{-\alpha}}{1-q^{\alpha -1}}\int\limits_{|y|_K< |x|_K}\left( |x|_K^{\alpha -1}-|y|_K^{\alpha -1}\right) u(|y|_K)\,dy,\quad \alpha \ne 1,
\end{equation}
and
\begin{equation}
(I^1 u)(|x|_K)=q^{-1}|x|_K u(|x|_K)+\frac{1-q}{q\log q}\int\limits_{|y|_K<|x|_K}\left( \log |x|_K-\log |y|_K\right) u(|y|_K)\,dy.
\end{equation}
\end{lem}

\medskip
{\it Proof}. It is sufficient to compute the integrals over the set $\{ y\in K:\ |y|_K=|x|_K\}$, and that is done using the integration formulas (2.3) and (2.5). $\qquad \blacksquare$

\medskip
It follows from Lemma 2 that the function $I^\alpha u$ is continuous if, for example, $u$ is bounded near the origin (see an estimate of the integral $I_{\alpha ,0}$ in the proof of Theorem 1 below). If $|u(|x|_K)|\le C|x|_K^{-\varepsilon}$, as $|x|_K\ge 1$, then $\left| \left( I^\alpha u\right) (|x|_K)\right| \le C|x|_K^{\alpha -\varepsilon}$, as $|x|_K\ge 1$. Here and below we denote by $C$ various (possibly different) positive constants.

It is easy to transform (3.10) and (3.11) further obtaining series involving $u(q^n)$.

Obviously, $D^\alpha I^\alpha =I$ on $\mathcal D(K)$, if $\alpha \ne 1$, and on $\Phi (K)$, if $\alpha =1$. Since by Lemma 1 and Lemma 2, the operators are defined in a straightforward sense for wider classes of functions, it is natural to look for conditions sufficient for this identity.

\medskip
\begin{lem}
Suppose that for some $m\in \mathbb Z$,
$$
\sum\limits_{k=-\infty}^m \max \left( q^k,q^{\alpha k}\right) \left| v(q^k)\right| <\infty ,\quad \sum\limits_{l=m}^\infty \left| v(q^l)\right| <\infty ,
$$
if $\alpha \ne 1$, and
$$
\sum\limits_{k=-\infty}^m |k|q^k \left| v(q^k)\right| <\infty ,\quad
\sum\limits_{l=m}^\infty l\left| v(q^l)\right| <\infty ,
$$
if $\alpha =1$. Then there exists $\left( D^\alpha I^\alpha v\right) (|x|_K)=v(|x|_K)$ for any $x\ne 0$.
\end{lem}

\medskip
The {\it proof} consists of tedious but quite elementary calculations based on the integration formulas (2.2)-(2.8). A relatively nontrivial tool is the sum formula for the arithmetic-geometric progression (formula 0.113 from \cite{GR}). $\qquad \blacksquare$

\medskip
Using Lemma 3, we can consider the simplest Cauchy problem
$$
D^\alpha u(|x|_K)=f(|x|_K),\quad u(0)=0,
$$
where $f$ is a continuous function, such that
$$
\sum\limits_{l=m}^\infty \left| f(q^l)\right| <\infty , \text{ if $\alpha \ne 1$, or } \sum\limits_{l=m}^\infty l\left| f(q^l)\right| <\infty ,\text{ if $\alpha =1$}.
$$
The unique strong solution is $u=I^\alpha f$; the uniqueness follows from the fact that the equality $D^\alpha u=0$ (in the sense of $\mathcal D'(K)$) implies the equality $u=\const$; see \cite{VVZ} or \cite{K}. Therefore on radial functions, the operators $D^\alpha$ and $I^\alpha$ behave like the Caputo-Dzhrbashyan fractional derivative and the Riemann-Liouville fractional integral of real analysis (see, for example, \cite{KST}). However the next example illustrates a different behavior of the ``fractional integral'' in the non-Archimedean case.
\begin{exam}
Let $f(|x|_K)\equiv 1$, $x\in K$. Then $\left(I^\alpha f\right) (|x|_K)\equiv 0$.
\end{exam}

\medskip
{\it Proof}. Let $|x|_K=q^n$. If $\alpha \ne 1$, then by (3.10), (2.2), and (2.6),
\begin{multline*}
\left(I^\alpha f\right) (|x|_K)=q^{-\alpha}|x|_K^\alpha +\frac{1-q^{-\alpha}}{1-q^{\alpha -1}}\int\limits_{|y|_K \le q^{n-1}}\left( |x|_K^{\alpha -1}-|y|_K^{\alpha -1}\right) \,dy\\
=q^{-\alpha}|x|_K^\alpha +\frac{1-q^{-\alpha}}{1-q^{\alpha -1}}\left[ q^{n-1}|x|_K^{\alpha -1}-\frac{1-q^{-1}}{1-q^{-\alpha }}q^{\alpha (n-1)}\right] \\
=q^{-\alpha}|x|_K^\alpha +\frac{1-q^{-\alpha}}{1-q^{\alpha -1}}|x|_K^\alpha \frac{q^{-1}-q^{-\alpha}}{1-q^{-\alpha }}=0.
\end{multline*}

If $\alpha =1$, then by (3.11), (2.4), and (2.6),
\begin{multline*}
\left(I^1 f\right) (|x|_K)=q^{-1}|x|_K+\frac{1-q}{q\log q}\int\limits_{|y|_K\le q^{n-1}}\left( \log |x|_K-\log |y|_K\right)\,dy\\
=q^{-1}|x|_K+\frac{1-q}{q\log q}\left[ q^{n-1}\log |x|_K-\left( n-1-\frac1{q-1}\right) q^{n-1}\log q\right] \\
=|x|_K\left( q^{-1}+\frac{1-q}{q\log q}\left( 1+\frac1{q-1}\right) q^{-1}\log q\right) =0. \qquad \blacksquare
\end{multline*}

\medskip
Of course, the above identities in the weaker sense of distributions from $\Phi'(K)$ are trivial, since the constant functions are identified with zero, $I^\alpha$ with $D^{-\alpha}$, and $D^\alpha D^{-\alpha }=I$.

On the other hand, the example shows that the condition of decay at infinity in Lemma 3 cannot be dropped.

\section{Fractional Differential Equations}

{\bf 4.1. The Cauchy problem and an integral equation.} In the class of radial functions $u=u(|x|_K)$, we consider the Cauchy problem
\begin{equation}
D^\alpha u+a(|x|_K)u=f(|x|_K),\quad x\in K,
\end{equation}
\begin{equation}
u(0)=0,
\end{equation}
where $a$ and $f$ are continuous functions, that is they have finite limits $a(0)$ and $f(0)$, as $x\to 0$.

Looking for a solution of the form $u=I^\alpha v$, where $v$ is a radial function, we obtain formally an integral equation
\begin{equation}
v(|x|_K)+a(|x|_K)\left( I^\alpha v\right) (|x|_K)=f(|x|_K),\quad x\in K.
\end{equation}

Let us study its solvability. Later we investigate, in what sense a solution of (4.3) corresponds to a solution of the Cauchy problem (4.1)-(4.2).

It follows from (4.3) that $v(0)=f(0)$. Suppose first that $\alpha \ne 1$. By Lemma 2, the equation (4.3) can be written in the form
\begin{multline*}
\left[ 1+q^{-\alpha}a(|x|_K)|x|_K^\alpha \right]v(|x|_K)\\
+c_\alpha a(|x|_K)\int\limits_{|y|_K<|x|_K}\left( |x|_K^{\alpha -1}-|y|_K^{\alpha -1}\right) v(|y|_K)\,dy=f(|x|_K),\quad x\ne 0,\tag{$4.3'$}
\end{multline*}
where $c_\alpha =\dfrac{1-q^{-\alpha}}{1-q^{\alpha -1}}$.

Since $a$ is continuous, there exists such $N\in \mathbb Z$ that
$$
q^{-\alpha}a(|x|_K)|x|_K^\alpha <1 \quad \text{for $|x|_K\le q^N$}.
$$
On the ball $B_N=\left\{ x\in K:\ |x|_K\le q^N\right\}$, the equation takes the form
\begin{equation}
v(|x|_K)+\int\limits_{|y|_K<|x|_K}k_\alpha (x,y)v(|y|_K)\,dy=F(|x|_K)
\end{equation}
where
$$
k_\alpha (x,y)=\left[ 1+q^{-\alpha}a(|x|_K)|x|_K^\alpha \right]^{-1}c_\alpha a(|x|_K)\left( |x|_K^{\alpha -1}-|y|_K^{\alpha -1}\right) ,\quad x\ne 0,
$$
$k_\alpha (0,y)=0$, $F(|x|_K)=\left[ 1+q^{-\alpha}a(|x|_K)|x|_K^\alpha \right]^{-1}f(|x|_K)$.

If we construct a solution of (4.4) on $B_N$, and if
\begin{equation}
a(|x|_K)\ne -q^{\alpha m}\quad \text{for any $x\in K$, $m\in \mathbb Z$},
\end{equation}
we will be able to construct a solution of (4.4), thus a solution of the equation (4.3), successively for all $x\in K$.

If $\alpha =1$, we use (3.11) and obtain in a similar way the equation (4.4) with
$$
k_1(x,y)=\frac{1-q}{q\log q}\left[ 1+q^{-1}a(|x|_K)|x|_K\right]^{-1}a(|x|_K)\left( \log |x|_K-\log |y|_K\right) ,\quad x\ne 0,
$$
$k_1(0,y)=0$, $F(|x|_K)=\left[ 1+q^{-1}a(|x|_K)|x|_K\right]^{-1}f(|x|_K)$.

It is convenient to extend $k_\alpha$ (including the case $\alpha =1$) by zero onto $B_N\times B_N$.

\medskip
\begin{teo}
For each $\alpha >0$, the integral equation (4.4) has a unique continuous solution on $B_N$.
\end{teo}

\medskip
{\it Proof}. Let us consider the integral operator $\mathcal K$ appearing in (4.4) as an operator on the Banach space $C(B_N)$ of complex-valued continuous functions on $B_N$. By the theory of integral operators developed in sufficient generality in \cite{E} (see Proposition 9.5.17), to prove that $\mathcal K$ is a compact operator, it suffices to check that, for any $x_0\in B_N$,
\begin{equation}
\lim\limits_{x\to x_0}\int\limits_{B_N}\left| k_\alpha (x,y)-k_\alpha (x_0,y)\right|\,dy=0.
\end{equation}

The relation (4.6) is obvious for $x_0\ne 0$, and also for $\alpha >1$. For $x_0=0$, we have $k_\alpha (0,y)=0$, and for $0<\alpha <1$, $|x|_K=q^n$, $n\ge N$, we get by (2.2) and (2.6) that
\begin{multline*}
\int\limits_{B_N}\left| k_\alpha (x,y)\right|\,dy=\const \cdot \int\limits_{|y|_K\le q^{n-1}}\left( |y|_K^{\alpha -1}-q^{n(\alpha -1)}\right)\,dy\\
=\const \cdot \left( \frac{1-q^{-1}}{1-q^{-\alpha}}q^{\alpha (n-1)}-q^{n(\alpha -1)}q^{n-1}\right)
=\const \cdot |x|_K^\alpha \to 0,
\end{multline*}
as $|x|_K\to 0$. For $\alpha =1$, we use (2.4) and (2.6) to obtain that
\begin{multline*}
\int\limits_{B_N}\left| k_1(x,y)\right|\,dy=\const \cdot \int\limits_{|y|_K\le q^{n-1}}\left( \log q^n-\log |y|_K\right)\,dy\\
=\const \cdot \left( nq^{n-1}\log q-(n-1-\frac1{q-1})q^{n-1}\log q\right)\\
=\const \cdot \frac{q}{q-1}q^{n-1}\log q
=\const \cdot \frac{\log q}{q-1}|x|_K\to 0,
\end{multline*}
as $|x|_K\to 0$.

Thus, $\mathcal K$ is compact, and by the Fredholm alternative (\cite{E}, 9.10.3), our theorem will be proved if we show that $\mathcal K$ has no nonzero eigenvalues.

Suppose that $\mathcal Kw=\lambda w$, $\lambda \ne 0$, for some $w\in C(B_N)$. We have $|w(y)|\le C$,
$$
|k_\alpha (x,y)|\le M\left| |x|_K^{\alpha -1}-|y|_K^{\alpha -1}\right|,
$$
if $\alpha \ne 1$, and
$$
|k_1(x,y)|\le M\left( \log |x|_K-\log |y|_K\right) ,
$$
if $\alpha =1$, $|y|_K<|x|_K$.

In subsequent iterations we will deal with the integrals
$$
I_{\alpha ,m}=\int\limits_{|y|_K<|x|_K}\left| |x|_K^{\alpha -1}-|y|_K^{\alpha -1}\right| |y|_K^{\alpha m}\,dy,\quad \alpha \ne 1,
$$
and
$$
I_{1,m}=\int\limits_{|y|_K<|x|_K}\left( \log |x|_K-\log |y|_K\right)|y|_K^m\,dy.
$$

If $\alpha >1$, we find denoting $|x|_K=q^n$ and using (2.2) that
$$
I_{\alpha ,m}=|x|_K^{\alpha -1}\int\limits_{|y|_K\le q^{n-1}}|y|_K^{\alpha m}\,dy-\int\limits_{|y|_K\le q^{n-1}}|y|_K^{\alpha (m+1)-1}\,dy=d_{\alpha ,m}|x|_K^{\alpha (m+1)}
$$
where, for all $m=0,1,2,\ldots$
\begin{multline*}
d_{\alpha ,m}=\frac{1-q^{-1}}{1-q^{-\alpha m-1}}q^{-\alpha m-1}-\frac{1-q^{-1}}{1-q^{-\alpha m-\alpha}}q^{-\alpha m-\alpha}\\
=(1-q^{-1})\frac{q^{\alpha -1}-1}{(1-q^{-\alpha m-1})(q^{\alpha m+\alpha}-1)}\le Aq^{-\alpha m}
\end{multline*}
for some $A>0$.
A similar result,
\begin{equation}
I_{\alpha ,m}=d_{\alpha ,m}|x|_K^{\alpha (m+1)},\quad d_{\alpha ,m}\le Aq^{-\alpha m},\quad m=0,1,2,\ldots
\end{equation}
is obtained for $0<\alpha <1$, so that (4.7) holds for any $\alpha \ne 1$. If $\alpha =1$, then the integral
$I_{1,m}$ is evaluated as follows. We have
\begin{multline*}
I_{1,m}=\sum\limits_{k=-\infty}^{n-1}\int\limits_{|y|_K=q^k}\left( \log |x|_K-\log |y|_K\right)|y|_K^m\,dy
=(1-\frac1q)\log q\sum\limits_{k=-\infty}^{n-1}(n-k)q^{k(m+1)}\\
=(1-\frac1q)\log q\sum\limits_{\nu =1}^\infty \nu q^{(n-\nu )(m+1)}=d_{1,m}|x|_K^{\alpha (m+1)}
\end{multline*}
where
$$
d_{1,m}=(1-\frac1q)\log q\sum\limits_{\nu =1}^\infty \nu q^{-\nu (m+1)}=(1-\frac1q)\log q\frac{q^{-m-1}}{(1-q^{-m-1})^2}\le Aq^{-m}
$$
(we have used the identity 0.231.2 from \cite{GR}). Thus, we have proved (4.7) also for $\alpha =1$.

Let us return to a function $w$ satisfying the relation $\mathcal Kw=\lambda w$, $\lambda \ne 0$. Using (4.7) (separately for $\alpha \ne 1$ and $\alpha =1$) and iterating we find by induction that
\begin{equation}
|w(x)|\le C\left( M|\lambda|^{-1}A\right)^m\left( \prod\limits_{j=0}^m q^{-\alpha j}\right) |x|_K^{\alpha m},\quad m=0,1,2,\ldots ,x\in B_N.
\end{equation}
Since $\prod\limits_{j=0}^m q^{-\alpha j}=q^{-\frac{\alpha}2 m(m+1)}$, it follows from (4.8) that $w(x)\equiv 0.\qquad \blacksquare$

\bigskip
{\bf 4.2. Strong solutions.} Below we assume that the inequality (4.5) is satisfied. Then, as we have mentioned, the solution $v$ of (4.4) is automatically extended in a unique way from $B_N$ onto $K$. The extended function $v$ satisfies the equation (4.3). Therefore the function $u=I^\alpha v$ satisfies (4.1) in the sense of distributions from $\Phi'$. The initial condition (4.2) is satisfied automatically.

Let us find additional conditions on $a$ and $f$, under which this construction gives a strong solution of the Cauchy problem (4.1)-(4.2). Note that, by Lemma 3 and Theorem 1, a strong solution is unique in the class of functions $u=I^\alpha v$ where $v$ is a continuous radial function, such that $\sum\limits_{l=m}^\infty \left| v(q^l)\right| <\infty$ for sum $m\in \mathbb Z$.

\medskip
\begin{teo}
Suppose that
\begin{equation}
\left| a(|x|_K)\right| \le C|x|_K^{-\alpha -\varepsilon},\quad \left| f(|x|_K)\right| \le C|x|_K^{-\varepsilon},\quad \varepsilon >0,C>0,
\end{equation}
as $|x|_K>1$. Then $u=I^\alpha v$ is a strong solution of the Cauchy problem (4.1)-(4.2).
\end{teo}

\medskip
{\it Proof}. Let $v(|x|_K)$ be the solution of the equation $(4.3')$ constructed above for all $x\in K$ (for $x=0$, the integral in the right-hand side is assumed equal to zero). For $|x|_K\le q^N$ the existence of a solution $v$ was obtained from the theory of compact operators; for larger values of $|x|_K$ we use successively the equation $(4.3')$ itself. Denote
$$
V_m=\sup\limits_{|x|_K\le q^m}\left| v(q^m)\right|.
$$
The sequence $\{V_m\}$ is non-decreasing.

As we assumed in Theorem 1 only the continuity of the coefficient $a$, we took $N$ in such a way that the neighborhood $B_N=\{ x:\ |x|_K\le q^N\}$ was sufficiently small. Here we assume (4.5), so that we can take any fixed integer $N$ and obtain a solution $v$ on $B_N$.

Consider the case where $\alpha \ne 1$. It follows from (4.5) and (4.9) that
$$
\left| \left[ 1+q^{-\alpha}a(|x|_K|x|_K^\alpha \right]^{-1}\right| \le H
$$
where $H>0$ does not depend on $x\in K$. If $m\ge N$, then we find from $(4.3')$ and the above estimate for $I_{\alpha ,0}$ that
\begin{equation}
\left| v(q^m)\right| \le c_\alpha d_{\alpha ,0}Ha(q^m)q^{\alpha m}V_{m-1}+H\left| f(g^m)\right|.
\end{equation}

Let us choose $m_1\ge N$ so big that
$$
c_\alpha d_{\alpha ,0}Ha(q^m)q^{\alpha m}\le \frac12,\quad H\left| f(g^m)\right|\le \frac12 V_{N-1},
$$
as $m\ge m_1$ (that is possible due to (4.9)). Then it follows from (4.10) that $V_m\le V_{m-1}$, as $m\ge m_1$, hence that the function $v$ is bounded on $K$.

Now we get from $(4.3')$ and the assumptions (4.9) that
\begin{equation}
\left| v(|x|_K)\right| \le C|x|_K^{-\varepsilon},\quad |x|_K\ge 1,
\end{equation}
$C>0$. A similar reasoning works for $\alpha =1$.

Taking into account the estimate (4.11) we find from Lemma 3 that $(D^\alpha I^\alpha v)(|x|_K)=v(|x|_K)$, $x\ne 0$. Therefore the function $u=I^\alpha v$ satisfies the equation (4.1) for all $x\ne 0$. Since $a,f$, and $u$ are continuous, the equation is satisfied in the strong sense. $\qquad \blacksquare$

\bigskip
{\bf 4.3. Generalizations.} Instead of (4.2), one can consider an inhomogeneous initial condition $u(0)=u_0$, $u_0\in \mathbb C$. Looking for a solution in the form $u=u_0+I^\alpha v$, $v=v(|x|_K)$, we obtain the integral equation
$$
v(|x|_K)+a(|x|_K)\left( I^\alpha v\right) (|x|_K)=f(|x|_K)-a(|x|_K)u_0,
$$
which can be studied under the same assumptions.

All the above results carry over to the case of a matrix-valued coefficient $a(|x|_K)$ and vector-valued solutions. In this case, to obtain a strong solution, it is sufficient to demand that the spectrum of each matrix $a(|x|_K)$, $x\in K$, does not intersect the set $\{-q^N,N\in \mathbb Z\}$.

\end{document}